\documentclass{llncs}
\usepackage{epsfig}

\usepackage{amsmath}
\usepackage{amsfonts}
\usepackage{amssymb}

\DeclareMathOperator{\sech}{sech}

\newcommand{\R}{{I\!\!R}}

\def\R{{\rm I}\! {\rm R}}

\renewcommand{\vec}[1]{\mbox{\boldmath $ #1 $}}

\newtheorem{algorithm}[theorem]{Algorithm}
\newtheorem{assumption}[theorem]{Assumption}

\begin{document}

\pagestyle{headings}

\title{Comparison of Splitting methods for Gross-Pitaevskii Equation}
\author{J\"urgen Geiser\inst{1} \and Amirbahador Nasari\inst{2}}
\institute{Ruhr University of Bochum, \\
Department of Electrical Engineering and Information Technology, \\
Universit\"atsstrasse 150, D-44801 Bochum, Germany \\
\email{juergen.geiser@ruhr-uni-bochum.de}
\and
Ruhr University of Bochum, \\
Department of Civil and Environmental Engineering, \\
Universit\"atsstrasse 150, D-44801 Bochum, Germany \\
\email{amirbahador.nasari@ruhr-uni-bochum.de}
}
\maketitle

\begin{abstract}

In this paper, we discuss the different splitting approaches
to solve the Gross-Pitaevskii equation numerically.
We consider conservative finite-difference schemes
and spectral methods for the spatial discretisation.
Further, we apply implicit or explicit time-integrators
and combine such schemes with different splitting approaches.
The numerical solutions are compared based on the conservation of
the $L_2$-norm with the analytical solutions.
The advantages of the splitting methods for large time-domains 
are presented in several numerical examples of different solitons
applications.

\end{abstract}



\section{Introduction}

Bose-Einstein condensate (BEC) nowadays is an actual modelling problem
for theoretical and also experimental studies, see \cite{dalfovo1999}.
The evolution equation of the Bose-Einstein condensate (BEC) order parameter for weakly interacting
bosons is done with the Gross-Pitaevskii equation, see \cite{abdullaev2005},
\cite{dauxois2006} and \cite{balak2011}.
The weakly interacting bosons supports dark solitons for repulsive interactions
and bright solitons for attractive interactions.
A solitary wave or soliton solution is a localised travelling wave solution,
that retain its size, shape and speed, when it moves. It does not spread or disperse,
see \cite{trofimov2009}. 
The modelling equation has two parts, a defocusing effect, which
is based on the dispersive term and a steeping effect, which is based on the nonlinear 
term. To obtain a equation balance of such a localised profile for the solution, we need
a special nonlinearity, see \cite{trofimov2009}. Also after a collision of two solitons, each wave is unscathed with its size, shape and speed, therefore,
we have a special collision property, see \cite{atre2006}. 
Therefore, the numerical methods should also have conservational behaviours 
to solve such a specialised balance of nonlinearity (steepness) and
diffusivity (smoothness) to obtain the sharp localised soliton solutions.

We are motivated to analyse such numerical methods, which allow
to conserve such behaviours, see \cite{trofimov2009} and \cite{jiang2013}.
Additionally, we apply different splitting approaches in combination
with finite difference schemes or spectral schemes to solve the 
Gross-Pitaevskii equation, see \cite{abdullaev2005}.

Numerically, two different ideas exist to solve the
GPE:
\begin{itemize}
\item Conservative finite difference schemes, which are nonlinear schemes and
need more computational amount, see \cite{trofimov2009}. Further, they can be constructed to
conserve the solution, the momentum and the energy. 
\item Splitting schemes, which decompose the different parts of the GPE and are simple to implement. But they have energy conservation and stability problems, see \cite{trofimov2009}
\end{itemize}

Based on the different ideas, we propose a combination of the splitting approaches and the uses of the conservation properties based on the conservative finite-difference schemes,
see \cite{trofimov2009}. Therefore, we could use the benefits of conservation approaches,
see \cite{trofimov2009} and the splitting approaches,
see \cite{geiser_2016} and \cite{geiser_2018_1} to stabilise and accelerate the solver processes.
Such a combination allows to reduce the time-consuming procedures
of the conservative FD schemes and stabilises the splitting approaches
based on the conservative approaches.

The paper is outlined as following. The model is introduced in Section \ref{modell}.
In Section \ref{methods}, we discuss the different numerical methods and present the
convergence analysis. The numerical experiments are done in Section \ref{numerics} and
the conclusion is presented in Section \ref{concl}.

\section{Mathematical Model}
\label{modell}

The modelling is based on many-body Hamiltonian for a system of $N$ interacting
particles (e.g., bosons) for the external field $V_{ext}$ and 
particle-particle interaction potential with $V(\vec{r} - \vec{r}')$:
\begin{eqnarray}
&& H = - \int \vec{u}^\dagger \left(\frac{\hbar^2}{2m} \; \nabla^2 - V_{ext}(\vec{r}) + \mu \right) \vec{u}  \; d\vec{r} + \nonumber \\
&& + \frac{1}{2} \int  \vec{u}^\dagger(\vec{r}) \vec{u}^\dagger(\vec{r}') \; V(\vec{r} - \vec{r}') \; \vec{u}(\vec{r})  \vec{u}(\vec{r}') \; d\vec{r}' d\vec{r}, 
\end{eqnarray}
where $\vec{u}$ is the particle (boson) field operator and we satisfy the commutation relation $[\vec{u}(\vec{r}), \vec{u}^\dagger(\vec{r}')] = \delta(\vec{r} - \vec{r}')$. Further, $V(\vec{r} - \vec{r}')$ is the two-body interaction
and $\mu$ is the chemical potential, see \cite{balak2011}. Then the time-evolution of the field operator $\vec{u}$ is given as:
\begin{eqnarray}
\label{boson_1}
 i \hbar \frac{\partial}{\partial t} \vec{u}(\vec{r}, t) && =  [\vec{u}, H] , \\
\label{boson_2}
 i \hbar \frac{\partial}{\partial t} \vec{u}(\vec{r}, t) && =  \left( - \frac{\hbar^2}{2m} \nabla^2 + V_{ext}(\vec{r}) - \mu + \right.\\
&& \left. +  \int \vec{u}^\dagger(\vec{r}', t) \; V(\vec{r} - \vec{r}') \; \vec{u}(\vec{r}',t) \; d\vec{r}'
\right)  \vec{u}(\vec{r}, t) . \nonumber
\end{eqnarray}

Further, the BEC order parameter, or called as condensate wave function, is given as $u = \langle \vec{u} \rangle$, where $\langle \vec{u} \rangle$ is the expectation value of the Bose operator. \\

We have two possibilities:
\begin{itemize}
\item $\langle \vec{u} \rangle = 0$, for $T > T_c$ and 
\item  $\langle \vec{u} \rangle \neq 0$, for $T < T_c$,
\end{itemize}
where $T_c$ is the Bose-Einstein condensation temperature.

In the following, we discuss the weakly interacting bosons.

\subsection{Weakly interacting Bosons}

We deal with the following Assumption \ref{assum_1}.
\begin{assumption}
\label{assum_1}

\begin{itemize}
\item We consider dilute gas, while we assume, that the range $r_0$ of the interatomic forces is much more smaller, than the distance between the atoms, means
$r_0 << d = n^{-1/3}$, where $n$ is the density of the atoms.
\item For $T < T_c$, we obtain small momenta, such that the scattering amplitude is independent of the energy. Therefore, one could replace it by a low-energy-value, which is determined by the solitary wave with scattering length $a$.
\item We replace the potential $V(\vec{r} - \vec{r}')$ with the effective soft potential $V_{eff}$, which has the same scattering properties, and it is defined as:
\begin{eqnarray}
 g =  \int V_{eff}(\vec{r}) \; d\vec{r} = \frac{4 \pi \; \hbar \; a}{m}
\end{eqnarray}
where $m$ is the atomic mass.
Further we replace $V(\vec{r} - \vec{r}') = g \; \delta(\vec{r} - \vec{r}')$.
\item We transform $\vec{u} \rightarrow \vec{u} \exp(i \; \mu \; t / \hbar)$.
\item The expectation value is given as $u = \langle \vec{u} \rangle$.
\end{itemize}

\end{assumption}

We apply the Assumption (\ref{assum_1}) to the evolution equation  of the interacting particle system (\ref{boson_2})
and obtain the Gross-Pitaevskii equation with the condensate order parameter $u$ for weakly interacting
bosons:
\begin{eqnarray}
&& i \hbar \frac{\partial u}{\partial t} =  \left( - \frac{\hbar^2}{2m} \; \nabla^2 + g |u(\vec{x}, t)|^2 \right) u , \; (\vec{x}, t) \in \R^3 \times [0, T] ,
\end{eqnarray}
where $g$ is the interaction term with the following characteristics:
\begin{itemize}
\item  $g > 0$ implies a repulsive interaction, where $a > 0$,
\item  $g < 0$ implies an attractive interaction, where $a < 0$.
\end{itemize}

The Gross-Pitaevskii equation is a nonlinear partial differential
equation with a cubic nonlinearity, means we deal with higher 
order nonlinearities, see also nonlinear Schr\"odinger equation \cite{tak2008}.

In the following, we concentrate on the one-dimensional Gross-Pitaevskii equation,
where we assume $\hbar = 1.0$ and the atomic mass $m=1$. We also deal with a external potential $V(x, t) \equiv 0$ and we deal with the
following form of the GPE:
\begin{eqnarray}
&& i \frac{\partial u}{\partial t} =  \left( - \frac{1}{2} \; \frac{\partial^2}{\partial x^2} + g |u(x, t)|^2 \right) u , \; (x, t) \in [-L, L] \times [0, T] , \\
&& u(x, t) = 0, \; x = \{-L, L\}, \; \mbox{and}  t \in [0, T] , \\
&& u(x, 0) = u_0(x), \; x \in [-L, L],
\end{eqnarray}
where the Hamiltonian operator is given as 
$H= \left( -\frac{1}{2} \; \frac{\partial^2}{\partial x^2} + g |u(x, t)|^2 \right)$.
Further, we assume $g = -1$, means we discuss attractive interactions.

\section{Numerical Methods}
\label{methods}

For the numerical methods, we deal with the two standard ideas
to approximate the GPE:
\begin{enumerate}
\item Splitting methods: The idea is to split the differential equations
into some simpler parts and solve each simpler differential equation
with fast PDE or ODE solvers. The results are summary approximated, e.g.,
via coupling the solution of the predecessor-solution as initial conditions 
of the successor-solution, or averaging the summarised results, see \cite{strang1968}
and \cite{geiser2011}.
The benefits are the fast solver methods and a simple numerical construction with the
simple implementation into a program-code, see \cite{strang1968} and \cite{mclach02}.
The drawback is that the methods are not long-time stable and they preserve only on invariant of the solution, see \cite{trofimov2009}. 
\item Conservative Finite Element Schemes: The idea is to design a finite difference scheme, which preserve the square of $L_2$-norm of the solution, the impulse functional and the energy functional. Based on such a construction of finite-difference approaches, e.g., a well-known conservative FD scheme is the 
semi-implicit Crank-Nicolson method, see \cite{sanz1986}, we conserve all the three invariants, see \cite{chang1999}, and we obtain stable and long-time behaviours of the solutions. The drawback of such schemes are the nonlinearity in the methods, e.g., we need additional nonlinear solvers, therefore the schemes are highly computational intensive comparing to fast splitting approaches, see \cite{trofimov2009}.
\end{enumerate}

We propose a mixture of the splitting approaches plus the application of the
conservation finite-difference schemes, while we
apply schemes for the GPE, which is given as:
\begin{eqnarray}
\label{gpe_1}
&& \frac{\partial u}{\partial t} = -i H u , \; x \in \Omega, \; t \in [0, 1] , \\
&& u(x,0) = u_0(x) , \; x \in \Omega ,  \\
\label{gpe_2}
&& u(x,t) = 0.0 , \; x \in \partial \Omega, \; t \in [0, 1] ,
\end{eqnarray}
with $ H u = \left(-\frac{1}{2} \frac{\partial^2}{\partial x^2} + g |u|^{2 \sigma} \right) u$, $\sigma = 1.0$ and we have applied Dirichlet boundary conditions.
Further, we apply $g=-1$, that means the attractive interaction case.

For an application of a single soliton, the exact solution is given as
\begin{eqnarray}
&& u(x,t) = A_0 \; \sech(\frac{|g|}{\sqrt{2}}(x - v_d \; t) A_0) \; \exp(i v_d (x - v_p \; t)/2) ,
\end{eqnarray}
where $A_0 = \sqrt{(v_d^2  - 2 v_p)/2 \; |g|}$, $v_d$ and $v_p$ are the speeds of the
density profile and phase profile, see the derivation of the exact solutions in 
\cite{balak2011}.

\begin{assumption}

We apply the absolute value as:
\begin{eqnarray}
&& | u(x,t) | = \sqrt{(\eta(x,t))^2 + (\xi(x,t))^2}.
\end{eqnarray}

Further we have the following complex relations:
\begin{eqnarray}
&& u(x,t) =  \eta(x, t) + i \xi(x, t) , \\
&& \exp(i \theta) = \cos(\theta) + i \sin(\theta), \\
&& \sech(\theta) = \frac{1}{\cosh(\theta)} .
\end{eqnarray}

\end{assumption}

\subsection{Conservation Laws of the GPE}

The GPE is given as in Equation (\ref{gpe_1})-(\ref{gpe_2}) and we have the following invariants:
\begin{itemize}
\item Mass conservation, which is given as the
square of $L_2$-norm of the solution
\begin{eqnarray}
{\cal N}(t) = \int_{\Omega} |u(x, t)|^2 \; dx ,
\end{eqnarray}
with ${\cal N}(t) = {\cal N}(0) = const$.

\item Impulse conservation, which is given as the
impulse functional of the solution
\begin{eqnarray}
{\cal P}(t) = \int_{\Omega} u^\dagger(x, t) (-i \frac{\partial}{\partial x}) u(x,t) \; dx ,
\end{eqnarray}
with ${\cal P}(t) = {\cal P}(0) = const$, with $u^\dagger$ is the conjugate of $u$.

\item Energy conservation, which is given as the
energy functional of the solution
\begin{eqnarray}
{\cal E}(t) = \frac{1}{2} \int_{\Omega} u^\dagger(x, t) \; H u(x,t) \; dx ,
\end{eqnarray}
with ${\cal E}(t) = {\cal E}(0) = const$, with $u^\dagger$ is the conjugate of $u$.

\end{itemize}

\begin{remark}
The conservation laws are proved for the GPE in the paper \cite{trofimov2009}.
Further, the conservation laws are also proved for the general Schr\"odinger equations in the paper \cite{chang1999}.
\end{remark}

In the following, we present a conservative finite difference scheme.

\subsection{Conservative finite difference schemes}

We apply the discretisation of the GPE (\ref{gpe_1})-(\ref{gpe_2})
with the following finite difference method, see also \cite{trofimov2009}:
\begin{eqnarray}
\label{cfds_1}
 i \frac{u^{n+1}_j - u^{n}_j}{\Delta t}&& =   - \frac{1}{2} \; \left( ( \frac{u_{j-1}^{n+1} - 2 u_{j}^{n+1} + u_{j+1}^{n+1}}{\Delta x^2} ) + ( \frac{u_{j-1}^{n} - 2 u_{j}^{n} + u_{j+1}^{n}}{\Delta x^2} ) \right) + \nonumber \\
&& + \frac{1}{2} \; g (|u_{j}^{n+1}|^2 + |u_{j}^{n}|^2)  \frac{u_{j}^{n+1} + u_{j}^{n}}{2} , \; j = 1, \ldots, M-1 , \\
&& u^0_j = u_0(x_j), \; j = 0, \ldots, M , \\
\label{cfds_2}
&& u^n_{0} = u_{M}^n = 0 , \; n = 0, 1, \ldots, N ,
\end{eqnarray}
where $M$ is the number of spatial grid points and $N$ is the number of time grid points.

Here, we have a conservative finite difference scheme, which has to be solved as a nonlinear equation
system with fixpoint or Newton's solvers, see \cite{trofimov2009}.
 
\begin{remark}
The conservative behaviour of the semi-implicit Crank-Nicolson is proved in \cite{trofimov2009}.
\end{remark}

\subsection{Asymptotic conservative finite difference schemes}

Here, we apply the idea of the conservative finite difference scheme
and reformulate the scheme into a splitting approach.

Therefore, we obtain asymptotic behaviours, while we have splitted the
full equations. Based on such a splitting approach, see \cite{geiser2011}, we have to apply additional
iterative steps to obtain the full coupled approximated conservative finite
difference scheme, see \cite{geiser2017_1}.

We reformulate the finite difference scheme (\ref{cfds_1})-(\ref{cfds_2}) in the
operator notation:
\begin{eqnarray}
\label{cfds_oper_1}
&&   U^{n+1} = U^n + i \; \frac{\Delta t}{2} \; \left( A_1 U^{n+1} + A_1 U^{n}  \right) + \nonumber \\
&& + i \; \frac{\Delta t}{2} \; \left( A_2(U^{n+1}) + A_2(U^n) \right) \frac{(U^{n+1} + U^n)}{2} , 
\end{eqnarray}
where the matrices are given as:
\begin{eqnarray}
&& A_1 = \frac{1}{2} \frac{1}{\Delta x^2 } \begin{bmatrix}
    -2 & 1 &   0 & 0 & \dots  & 0 \\
     1 & -2 &  1 & 0 & \dots  & 0 \\
     0 &  1 & -2 & 1 & \dots  & 0 \\
    \vdots & \vdots & \vdots & \vdots & \ddots & \vdots \\
    0 & 0 & 0 & 0 & \dots  & -2
\end{bmatrix} \in \R^{M-1 \times M-1} ,  \\
\label{cfds_oper_2}
&& A_2(U^n) = - g \; I \; abs(U^n)^2 \in \R^{M-1 \times M-1} ,
\end{eqnarray}
where with $U^n = (u_1^n, \ldots, u_{M-1}^n)^t$ is the vector at the grid points $u_j^n = u^n(x_j)$ for $j = 1, \ldots, M-1$. 
Further $I \in \R^{M-1 \times M-1}$ is the identity matrix and  $abs(U^n)^2 = (|u_1^n|^2, \ldots, |u_{M-1}^n|^2)^t \in R^{M-1}$ is a vector.

Further, the time-steps are given as $\Delta t = t^{n+1} - t^n$, with $n = 0, \ldots, N-1$ and $t^0 = 0$ and $i$ is the imaginary number.

We apply the following asymptotic approximation, based on the
Picards-fixpoint scheme, we reformulate the operator scheme (\ref{cfds_oper_1})-(\ref{cfds_oper_2}) 
as following:
\begin{eqnarray}
&&   U_k^{n+1} = U^n + i \; \frac{\Delta t}{2} \; \left( A_1 U_k^{n+1} + A_1 U^{n}  \right) + \nonumber \\
\label{asymp_cfds}
&& + i \; \frac{\Delta t}{2} \; \left( A_2(U_{k-1}^{n+1}) + A_2(U^n) \right) \frac{(U_{k-1}^{n+1} + U^n)}{2} , 
\end{eqnarray}
where $k = 1, \ldots, K$ is the iteration index and we have $U_{0}^{n+1} = U^n$ as the initialisation of the iteration, while we
have the stopping criterion $|| U_k^{n+1} -  U_{k-1}^{n+1} || \le err$ and $err$ is an error-bound, e.g., $err = 10^{-5}$, or we stop at $k = K$,
while $K$ is a fixed integer, e.g., $K = 5$.

We reformulate in a scaled $\frac{1}{2} AB$ and $\frac{1}{2} BA$ splitting approach.
Here, we obtain a first order splitting approach for both splitting approaches, see \cite{geiser2011}, see the Algorithm \ref{algo_aba_iter}.

\begin{algorithm}
\label{algo_aba_iter}
 We apply the time-steps $n= 1, \ldots, N-1$, where $N$ are the number of the time-steps.
The initialisation is $U^0 = U(0)$ and we start with $n=1$.

\begin{enumerate}
\item $\frac{1}{2} AB$
\begin{eqnarray}
&&   \tilde{U}_k^{n+1} = U^n + i \; \frac{\Delta t}{2} \; A_1 \tilde{U}_k^{n+1} + \nonumber \\
&& + i \; \frac{\Delta t}{2} \; \left( A_2(U_{k-1}^{n+1}) + A_2(U^n) \right) \frac{U_{k-1}^{n+1}}{2} , 
\end{eqnarray}
where the starting condition at $k=1$ is $U^{n+1}_0 = U^n$.
\item $\frac{1}{2} BA$
\begin{eqnarray}
&&   \hat{U}_k^{n+1} = \hat{U}^n + i \; \frac{\Delta t}{2} \; A_1 \hat{U}^{n} + \nonumber \\
&& + i \; \frac{\Delta t}{2} \; \left( A_2(U_{k-1}^{n+1}) + A_2(\hat{U}^n) \right) \frac{\hat{U}^n}{2} , 
\end{eqnarray}
where the starting condition at $k=1$ is $U^{n+1}_0 = U^n$, further we have $\hat{U}^n = \tilde{U}_k^{n+1}$.
The solution is given as $U_k^{n+1} = \hat{U}_k^{n+1}$.\\

If $k = K$ or $||U_k^{n+1} - U_{k-1}^{n+1} || \le err$, we are done and goto step 3., \\
else we go to the next iterative-step and we apply $k = k+1$ and goto step 1.

\item If $n+1 = N$, we are done, \\
else go to the next time-step and we apply $n = n+1$ and goto step 1.

\end{enumerate}

\end{algorithm}

We solve the two $B$-steps exactly and reformulate the asymptotic conservative finite difference scheme (\ref{asymp_cfds}) with respect to the
splitting approach, we call it the A-B-A(semiCN) splitting approach, see the Algorithm \ref{aba-semicn}.

Here the $A$ operator is the linear term with the FD scheme discretised,
while the $B$ operator is the nonlinear term and is exactly solved.
We apply an additional iterative procedure to approach the semi-implicit CN method.
\begin{algorithm}
\label{aba-semicn}

\begin{eqnarray}
&&   U_1^{n+1} = (I - i \Delta t/2 \; A_1)^{-1}  \; U^n , \mbox{with timestep} \; \Delta t/2 \; \mbox{(implicit Euler)} , \\
&&   U_2^{n+1} = \exp(- i \; g \; A_2 \Delta t) \; U_1^{n+1} , \mbox{with timestep} \; \Delta t \; \mbox{(spectral method)}, \\
&&   U^{n+1}_i = (I + i \Delta t/2 \; A_1)  \; U_2^{n+1} , \mbox{with timestep} \; \Delta t/2 \; \mbox{(explicit Euler)} ,
\end{eqnarray}
where
\begin{eqnarray}
&& A_1(t,x) = \; \frac{1}{2} \frac{1}{\Delta x^2 } \begin{bmatrix}
    -2 & 1 &   0 & 0 & \dots  & 0 \\
     1 & -2 &  1 & 0 & \dots  & 0 \\
     0 &  1 & -2 & 1 & \dots  & 0 \\
    \vdots & \vdots & \vdots & \vdots & \ddots & \vdots \\
    0 & 0 & 0 & 0 & \dots  & -2
\end{bmatrix} \in \R^{M-1 \times M-1} ,  \\
&& \hspace{-0.5cm} A_2(t,x, U^n, U_{i-1}^{n+1}) = \; I \frac{1}{2} \; \left( abs((U^n))^2 +  abs((U_{i-1}^{n+1}))^2  \right)  \in \R^{M-1 \times M-1} ,
\end{eqnarray}
where with spatial vector $x = (x_1, \ldots, x_{M-1})^t$ and $M$ are the number of spatial points.
Further $U^n = (u_1^n, \ldots, u_{M-1}^n)^t$ is the vector at the grid points $u_j^n = u^n(x_j)$ for $j = 1, \ldots, M-1$. 

The starting condition for $ U^{n+1}_0 =  U^n$.

\end{algorithm}

\begin{remark}
We reformulated the semi-CN scheme into an ABA-splitting approach, while the 
reformulation has also second order terms, we have at least for such
an approximation, only a first order scheme, see \cite{geiser2011}.
\end{remark}

\subsection{Standard Finite Difference Methods and Standard Splitting Approaches}

In the following, we discuss the different standard finite difference method
and standard Splitting approaches, which are related to the finite
difference schemes for the Gross-Pitaevskii equation.

\subsubsection{Splitting methods with finite difference schemes}

We apply the semi-discretisation of the diffusion
operator with a finite difference scheme (second order),
where we deal with $M$ discrete spatial points.

Further, We employ the following transformation and change of variables
with $u = \eta + i \xi \in (\R^M + i \R^M)$ and obtain:
\begin{eqnarray}
&&   U^{n+1} = U^n + i \; \Delta t \; A(t,x, U^n) U^n \\
&& A(t,x, U^n) = A_1(t, x) + A_2(t,x, U^n) , \\
&& A_1(t,x) = \frac{1}{2} \frac{1}{\Delta x^2 } \begin{bmatrix}
    -2 & 1 &   0 & 0 & \dots  & 0 \\
     1 & -2 &  1 & 0 & \dots  & 0 \\
     0 &  1 & -2 & 1 & \dots  & 0 \\
    \vdots & \vdots & \vdots & \vdots & \ddots & \vdots \\
    0 & 0 & 0 & 0 & \dots  & -2
\end{bmatrix} \in \R^{M-1 \times M-1} ,  \\
&& A_2(t,x, U^n) = \epsilon \; I \; abs(U^n)^2 \in \R^{M-1 \times M-1} ,
\end{eqnarray}
where with $U^n = (u_1^n, \ldots, u_{M-1}^n)^t$ is the vector at the grid points $u_j^n = u^n(x_j)$ for $j = 1, \ldots, M-1$. 

Further, the time-steps are given as $\Delta t = t^{n+1} - t^n$, with $n = 0, \ldots, N-1$ and $t^0 = 0$ and $i$ is the imaginary number.

\begin{itemize}

\item Implicit Euler method:
\begin{eqnarray}
  &&  U^{n+1} = (I - i \; \Delta t \; A(t,x, U^n))^{-1} U^n ,
\end{eqnarray}
where, we start with $U^0$.

\item CN-method:
\begin{eqnarray}
  &&  U^{n+1} = (I - i \; \Delta t/2 \; A(t,x, U^n))^{-1} (I - i \; \Delta t/2 \; A(t,x, U^n))^{-1} U^n ,
\end{eqnarray}
where, we start with $U^0$.
  
\item A--B splitting, where we deal with implicit for the diffusion and explicit time discretisation for the nonlinear term:
\begin{eqnarray}
  &&   U^{n+1} = U^n + i \Delta t \; ( A_1(t, x) U^{n+1} + A_2(t,x, U^n) ) U^n , \\
  &&  U^{n+1} = ( I - i \Delta t \; A_1(t, x) )^{-1} (I + i \Delta t \; A_2(t,x, U^n) )  U^n  ,
\end{eqnarray}
where we start with $U^0$.

\item A--B splitting, where we deal with explicit for the diffusion and explicit time discretisation for the nonlinear term:
\begin{eqnarray}
  &&   U^{n+1} = U^n + i \Delta t \; ( A_1(t, x) U^{n} + A_2(t,x, U^n) ) U^n , \\
  &&  U^{n+1} = ( I + i  \Delta t \; A_1(t, x) + i \Delta t \; A_2(t,x, U^n) )  U^n  ,
\end{eqnarray}
where we start with $U^0$.

\end{itemize}

\subsection{Standard Spectral Methods and Combinations with Splitting and Finite Difference schemes}

In the following, we present spectral and mixed schemes, 
combing spectral and finite difference schemes with splitting approaches.

The spectral methods applied the Fourier transformation
or Fourier spectral method, see \cite{brigham1973}.
The spectral methods can be applied to the linear part (spatial derivation)
and nonlinear part (interaction or potential) of the GPE, see \cite{trofimov2009}.

In the following, we apply the different splitting approaches with respect to
the spectral methods.

\subsubsection{Time-spitting spectral method}

We apply the spectral method in $t \in [t^n, t^{n+1}]$

We have two parts of the equation:
\begin{itemize}
\item Linear part:
\begin{eqnarray}
&& \frac{\partial u}{\partial t} = i \; \frac{1}{2} \frac{\partial^2 u}{\partial x^2} , (x, t) \in [-L, L] \times [0, T], \\
&&  u(x, t) = 0, \; x \in \{-L, L\}, \; t \in [0, T] ,
\end{eqnarray}
where we start to apply the Fourier transform for the input $u^n$ and
obtain:
\begin{eqnarray}
&& \hat{u}^n = \sum_{j = -M+1}^{M-1}  u_j^n \exp(- i \; \mu_l (x_j - L)) , \; l = - \frac{M}{2}, \ldots, \frac{M}{2}-1 , \\
&& \mu_l = \frac{\pi \;l}{L}, \; l = - \frac{M}{2}, \ldots, \frac{M}{2}-1 . 
\end{eqnarray}
We apply the Fourier transform to the linear term and obtain the result in the Fourier transformed space
and the inverse Fourier transform and obtain the result:
\begin{eqnarray}
&& u^{n+1} = \frac{1}{M} \sum_{l = - M/2}^{M/2 - 1} \exp(- i \; \mu^2_l \frac{\Delta t}{2}) \; \hat{u}_l^n \; \exp(i \; \mu_l (x_j - L)) 
\end{eqnarray}
\item Nonlinear part:
\begin{eqnarray}
&& \frac{\partial u}{\partial t} = - i \; g |u|^2 \; x \in \{-L, L\}, \; t \in [0, T] ,
\end{eqnarray}
where we obtain an analytical solution, which is given as:
\begin{eqnarray}
&& u^{n+1} = \exp(- i \; g |u^n|^2 \Delta t) \; u^n ,
\end{eqnarray}
where $\Delta t = t^{n+1} - t^n$.
\end{itemize}

The algorithm for the splitting approach is given as:

\begin{algorithm}

We apply the Time-splitting spectral method as following:
\begin{eqnarray}
&&   U_1^{n+1/2} = \exp(- i \; g |u^n|^2 \Delta t/2) \; U^n , \mbox{with timestep} \; \Delta t/2 , \\
&&   U_2^{n+1} = \frac{1}{M} \sum_{l = - M/2}^{M/2 - 1} \exp(- i \; \mu^2_l \frac{\Delta t}{2}) \; \hat{U}_{1, l}^{n+1/2} \; \exp(i \; \mu_l (x_j - L)) , \nonumber \\
&& \mbox{with timestep} \; \Delta t , \\
&&   U^{n+1} =  \exp(- i \; g |u^n|^2 \Delta t/2) \; U_2^{n+1} , \mbox{with timestep} \; \Delta t/2 ,
\end{eqnarray}
where $\hat{U}_{1}^{n+1/2} = \sum_{j = -M+1}^{M-1}  U_{1,j}^{n+1/2} \exp(- i \; \mu_l (x_j - L)) , \; l = - \frac{M}{2}, \ldots, \frac{M}{2}-1$ and $\mu_l = \frac{\pi \;l}{L}, \; l = - \frac{M}{2}, \ldots, \frac{M}{2}-1$.

Then, we start again with $U^{n+1}$ in step A.

\end{algorithm}

\subsubsection{AB Splitting Methods with finite difference and spectral schemes}

We deal with the different AB-splitting methods:
\begin{itemize}

\item 1.) TSSP Method: A and B are in the spectral version 

\item 2.) A-B splitting: A operator is the nonlinear term with the spectral method for the reaction \\
                         B operator is the linear term and is in the FD scheme
\item 3.) A-B splitting:   A operator is the nonlinear term with the FD scheme \\
                           B operator is the linear term in spectral method
\item 4.) A-B splitting:   A operator is the nonlinear term with the FD scheme
                            B operator is the linear term is in FD scheme 
\end{itemize}

\begin{itemize}

\item 1.) TSSP Method: A and B are in the spectral version 

\begin{algorithm}

We apply the Time-splitting spectral method as following:
\begin{eqnarray}
&&   U_1^{n+1} = \exp(- i \; g |u^n|^2 \Delta t) \; U^n , \mbox{with timestep} \; \Delta t ,\\
&&   U^{n+1} = \frac{1}{M} \sum_{l = - M/2}^{M/2 - 1} \exp(- i \; \mu^2_l \frac{\Delta t}{2}) \; \hat{U}_{1, l}^{n+1/2} \; \exp(i \; \mu_l (x_j - L)) , \\
&& \mbox{with timestep} \; \Delta t , \nonumber
\end{eqnarray}
where $\hat{U}_{1}^{n+1} = \sum_{j = -M+1}^{M-1}  U_{1,j}^{n+1} \exp(- i \; \mu_l (x_j - L)) , \; l = - \frac{M}{2}, \ldots, \frac{M}{2}-1$ and $\mu_l = \frac{\pi \;l}{L}, \; l = - \frac{M}{2}, \ldots, \frac{M}{2}-1$.
Then, we start again with $U^{n+1}$ in step A.

\end{algorithm}

\item 2.) A-B splitting: A operator is the nonlinear term with the spectral method for the reaction \\
                         B operator is the linear term and is in the FD scheme.

\begin{algorithm}
We apply the combined FD and spectral method as:
\begin{eqnarray}
&&   U_1^{n+1} = \exp(- i \; g |u^n|^2 \Delta t) \; U^n , \mbox{with timestep} \; \Delta t , \\
&&   U^{n+1} = ( I +  \; A_1(t,x) ) U^{n+1}_1 , \; \mbox{with timestep} \; \Delta t,
\end{eqnarray}
where %
\begin{eqnarray}
&& A_1(t,x) = i \; \frac{1}{2} \frac{\Delta t}{\Delta x^2 } \begin{bmatrix}
    -2 & 1 &   0 & 0 & \dots  & 0 \\
     1 & -2 &  1 & 0 & \dots  & 0 \\
     0 &  1 & -2 & 1 & \dots  & 0 \\
    \vdots & \vdots & \vdots & \vdots & \ddots & \vdots \\
    0 & 0 & 0 & 0 & \dots  & -2
\end{bmatrix} \in \R^{M-1 \times M-1} .
\end{eqnarray}

Then, we start again with $U^{n+1}$ in step A.

\end{algorithm}

\item 3.) A-B splitting:   A operator is the nonlinear term with the FD scheme \\
                           B operator is the linear term in spectral method

\begin{algorithm}

We apply the Time-splitting spectral method as following:
\begin{eqnarray}
&&   U_1^{n+1} =  U^n + (- i \; g A_2 \Delta t) \; U^n , \mbox{with timestep} \; \Delta t , \\
&&   U^{n+1} = \frac{1}{M} \sum_{l = - M/2}^{M/2 - 1} \exp(- i \; \mu^2_l \frac{\Delta t}{2}) \; \hat{U}_{1, l}^{n+1/2} \; \exp(i \; \mu_l (x_j - L)) , \\
&& \mbox{with timestep} \; \Delta t , \nonumber
\end{eqnarray}
where $\hat{U}_{1}^{n+1} = \sum_{j = -M+1}^{M-1}  U_{1,j}^{n+1} \exp(- i \; \mu_l (x_j - L)) , \; l = - \frac{M}{2}, \ldots, \frac{M}{2}-1$ and $\mu_l = \frac{\pi \;l}{L}, \; l = - \frac{M}{2}, \ldots, \frac{M}{2}-1$
and
\begin{eqnarray}
&& \hspace{-2cm} A_2(t,x, U) = \;  \begin{bmatrix}
    f(\eta_1, \xi_1, t^{n}, x_1) & 0 &   0 & 0 & \dots  & 0 \\
     0 & f(\eta_2, \xi_2, t^{n}, x_2) &  0 & 0 & \dots  & 0 \\
     0 &  0 & f(\eta_3, \xi_3, t^{n}, x_3) & 0 & \dots  & 0 \\
    \vdots & \vdots & \vdots & \vdots & \ddots & \vdots \\
    0 & 0 & 0 & 0 & \dots  & f(\eta_{M-1}, \xi_{M-1}, t^{n}, x_{M-1})
\end{bmatrix} \nonumber \\
&& \in \R^{M-1 \times M-1} ,
\end{eqnarray}
where $f(\eta_j, \xi_j, t^{n}, x_j) = (\eta(t^{n},x_j))^2 + (\xi(t^{n},x_j))^2$
for $j = 1, \ldots, M-1$ with the spatial vector $x = (x_1, \ldots, x_{M-1})$ and $M$ are the number of spatial points.
Further $U = (u_1, \ldots, u_{M-1})^t$ is the vector at the grid points $u_j = u(x_j)$ for $j = 1, \ldots, M-1$. 

Then, we start again with $U^{n+1}$ in step A.

\end{algorithm}

\item 4.) A-B splitting:   A operator is the nonlinear term with the FD scheme
                            B operator is the linear term is in FD scheme

\begin{algorithm}

We apply the splitting approach with the FD schemes as:
\begin{eqnarray}
&&   U_1^{n+1} = U^n + (- i \; g A_2 \Delta t) \; U^n , \mbox{with timestep} \; \Delta t , \\
&&   U^{n+1} = U^n + (i A_1 \Delta t) \; U_1^{n+1} , \mbox{with timestep} \; \Delta t ,
\end{eqnarray}
where
\begin{eqnarray}
&& A_1(t,x) = \; \frac{1}{2} \frac{1}{\Delta x^2 } \begin{bmatrix}
    -2 & 1 &   0 & 0 & \dots  & 0 \\
     1 & -2 &  1 & 0 & \dots  & 0 \\
     0 &  1 & -2 & 1 & \dots  & 0 \\
    \vdots & \vdots & \vdots & \vdots & \ddots & \vdots \\
    0 & 0 & 0 & 0 & \dots  & -2
\end{bmatrix} \in \R^{M-1 \times M-1} ,  \\
&& A_2(t,x, U^n) = \; I \; abs((U^n))^2  \in \R^{M-1 \times M-1} ,
\end{eqnarray}
where with spatial vector $x = (x_1, \ldots, x_{M-1})$ and $M$ are the number of spatial points.
Further $U^n = (u_1^n, \ldots, u_{M-1}^n)^t$ is the vector at the grid points $u_j^n = u^n(x_j)$ for $j = 1, \ldots, M-1$. 
\end{algorithm}

\item 5.) A-B-A(CN) splitting:   A operator is the linear term with the FD scheme
                            B operator is the nonlinear term is in spectral method

\begin{algorithm}

We apply the ABA-splitting approach with FD schemes and spectral schemes as:
\begin{eqnarray}
&&   U_1^{n+1} = (I - i \Delta t/2 \; A_1)^{-1}  \; U^n , \mbox{with timestep} \; \Delta t/2 \mbox{(implicit Euler)} , \\
&&   U_2^{n+1} = \exp(- i \; g \; A_2 \Delta t) \; U_1^{n+1} , \mbox{with timestep} \; \Delta t \mbox{(spectral method)}, \\
&&   U^{n+1} = (I + i \Delta t/2 \; A_1)  \; U_2^{n+1} , \mbox{with timestep} \; \Delta t/2 \mbox{(explicit Euler)} ,
\end{eqnarray}
where
\begin{eqnarray}
&& A_1(t,x) = \; \frac{1}{2} \frac{1}{\Delta x^2 } \begin{bmatrix}
    -2 & 1 &   0 & 0 & \dots  & 0 \\
     1 & -2 &  1 & 0 & \dots  & 0 \\
     0 &  1 & -2 & 1 & \dots  & 0 \\
    \vdots & \vdots & \vdots & \vdots & \ddots & \vdots \\
    0 & 0 & 0 & 0 & \dots  & -2
\end{bmatrix} \in \R^{M-1 \times M-1} ,  \\
&& A_2(t,x, U^n) = \; I \; abs((U^n))^2  \in \R^{M-1 \times M-1} ,
\end{eqnarray}
where with spatial vector $x = (x_1, \ldots, x_{M-1})$ and $M$ are the number of spatial points.
Further $U^n = (u_1^n, \ldots, u_{M-1}^n)^t$ is the vector at the grid points $u_j^n = u^n(x_j)$ for $j = 1, \ldots, M-1$. 
\end{algorithm}

\end{itemize}

\section{Numerical experiments}
\label{numerics}

For the numerical experiments, we test two models:
\begin{itemize}
\item Single soliton with exact solution as corresponding solution.
\item Collision of two solitons with numerically fine solution as corresponding solution.
\end{itemize}

For the errors, we apply the $L_2$-norm and use:

\begin{eqnarray}
err_{L_2, num, \Delta x, \Delta t} && = \left( \int_{[0, T]} \int_{\Omega} || u_{exact}(x, t) - u_{num}(x, t) ||^2 dx \; dt \; \right) = \nonumber \\
&& = \left( \Delta t \; \Delta x \; \sum_{n=1}^N \sum_{i=1}^{M} || u_{exact}(x_i, t^n) - u_{num}(x_i, t^n) ||^2  \; \right) , 
\end{eqnarray}
where $ || u_{exact}(x_i, t^n) - u_{num}(x_i, t^n) || = \mbox{\rm abs}( u_{exact}(x_i, t^n) - u_{num}(x_i, t^n))$.

We apply a convergence-tableau based on the different spatial- and time-steps,
means we apply $16 \Delta t, \ldots, \Delta t/8$ and  $\Delta x, \ldots, \Delta x/8$ with the underlying errors.

In the following, we apply different numerical experiments to
validate our numerical method.

\subsection{First example: GPE with one soliton}

We consider the GPE in order to apply for the numerical schemes in a suitable
rewriting:
\begin{eqnarray}
&& \frac{\partial u}{\partial t} = -i H u , \; x \in \Omega, \; t \in [0, 1] , \\
&& u(x,0) = \sech(\frac{1}{\sqrt{2}}(x -25)) \exp(i \frac{x}{20}) , \; x \in \Omega ,  \\
&& u(x,t) = 0.0 , \; x \in \partial \Omega, \; t \in [0, 1] ,
\end{eqnarray}
with $ H u = \left(-\frac{1}{2} \frac{\partial^2}{\partial x^2} + g |u|^{2 \sigma} \right) u$, $\sigma = 1.0$ and we have applied Dirichlet boundary conditions.

We applied for the analytical solution $g= -1$, $v_d = \frac{1}{10}$ and $v_p = -\frac{199}{200}$
and the analytical solution is given as:
\begin{eqnarray}
\hspace{-1cm} u(x,t) = \sech(\frac{1}{\sqrt{2}}(x - \frac{t}{10} - 25)) \exp(i (\frac{x}{20} - \frac{199}{400} t)) , \; (x, t) \in [-L, L] \times [0, T] .
\end{eqnarray}

We deal with the following methods:
\begin{itemize}
\item implicit Euler method (all operators are done with the implicit method),
\item Crank-Nicolson scheme (all operators are done with the CN method),
\item AB-splitting:
\begin{itemize}
\item linear operator is done with the Spectral method and nonlinear operator is done with the spectral method,
\item linear operator is done with the FD method and nonlinear operator is done with the spectral method,
\item linear operator is done with the Spectral method and nonlinear operator is done with the FD method,
\item linear operator is done with the FD method and nonlinear operator is done with the FD method.
\end{itemize}
\item ABA-splitting:
\begin{itemize}
\item linear operator is done with the Spectral method and nonlinear operator is done with the spectral method.
\end{itemize}
\item ABA-CN and ABA-iCN:
\begin{itemize}
\item linear operator is done with the finite difference method,
while the nonlinear operator is done with the spectral method.
\item for the iterative scheme, we apply different iterative steps.
\end{itemize}

\end{itemize}

The convergence-tableaus of the different numerical methods are given in the
Tables \ref{table_conv_1}-\ref{table_conv_10}.

\begin{table}
\centering
\begin{tabular}{ |c|c|c|c| } 
 \hline
  & $\Delta x/4$ & $\Delta x/8$ & $\Delta x/16$  \\ 
\hline
$4\Delta t $&1.5474e-06&3.9236e-11&1.7211e-13\\
\hline
$8\Delta t $&3.667e-05&9.2733e-10&4.6629e-12\\
\hline
$16\Delta t $&7.334e-05&1.8547e-09&9.3258e-12\\
\hline
\end{tabular}
\caption{\label{table_conv_1} Convergence tableau for the method implicit Euler.}
\end{table}

\begin{table}
\centering
\begin{tabular}{ |c|c|c|c| } 
 \hline
  & $\Delta x/4$ & $\Delta x/8$ & $\Delta x/16$  \\ 
\hline
$4\Delta t $&1.312e-05&3.3337e-10&1.6673e-12\\
\hline
$8\Delta t $&3.667e-05&9.2733e-10&4.6629e-12\\
\hline
$16\Delta t $&7.334e-05&1.8547e-09&9.3258e-12\\
\hline
\end{tabular}
\caption{\label{table_conv_2} Convergence tableau for the method Crank-Nicolson.}
\end{table}

\begin{table}
\centering
\begin{tabular}{ |c|c|c|c| } 
 \hline
 & $\Delta x/4$ & $\Delta x/8$ & $\Delta x/16$  \\ 
\hline
$4\Delta t $&1.583e-05&4.2155e-10&2.2106e-12\\
\hline
$8\Delta t $&3.667e-05&9.2733e-10&4.6629e-12\\
\hline
$16\Delta t $&7.334e-05&1.8547e-09&9.3258e-12\\
\hline
\end{tabular}
\caption{\label{table_conv_3} Convergence tableau for the method AB-splitting: A and B operators are spectral.}
\end{table}

\begin{table}
\centering
\begin{tabular}{ |c|c|c|c| } 
 \hline
 & $\Delta x/4$ & $\Delta x/8$ & $\Delta x/16$  \\ 
\hline
$4\Delta t $&1.312e-05&3.3337e-10&1.6673e-12\\
\hline
$8\Delta t $&3.667e-05&9.2733e-10&4.6629e-12\\
\hline
$16\Delta t $&7.334e-05&1.8547e-09&9.3258e-12\\
\hline
\end{tabular}
\caption{\label{table_conv_4}Convergence tableau for the method  AB-splitting: A spectral , B FD.}
\end{table}

\begin{table}
\centering
\begin{tabular}{ |c|c|c|c| } 
 \hline
  & $\Delta x/4$ & $\Delta x/8$ & $\Delta x/16$  \\ 
\hline
$4\Delta t $&1.583e-05&4.2155e-10&2.2106e-12\\
\hline
$8\Delta t $&3.667e-05&9.2733e-10&4.6629e-12\\
\hline
$16\Delta t $&7.334e-05&1.8547e-09&9.3258e-12\\
\hline
\end{tabular}
\caption{\label{table_conv_5} Convergence tableau for the method  AB-splitting: A FD , B spectral.}
\end{table}

\begin{table}
\centering
\begin{tabular}{ |c|c|c|c| } 
 \hline
 &$\Delta x/4$ & $\Delta x/8$ & $\Delta x/16$  \\ 
\hline
$4\Delta t $&1.312e-05&3.3337e-10&1.6673e-12\\
\hline
$8\Delta t $&3.667e-05&9.2733e-10&4.6629e-12\\
\hline
$16\Delta t $&7.334e-05&1.8547e-09&9.3258e-12\\
\hline
\end{tabular}
\caption{\label{table_conv_6}Convergence tableau for the method  AB-splitting: A FD , B FD. }
\end{table}

\begin{table}
\centering
\begin{tabular}{ |c|c|c|c| } 
 \hline
  & $\Delta x/4$ & $\Delta x/8$ & $\Delta x/16$  \\ 
\hline
$4\Delta t $&1.583e-05&4.2155e-10&2.2106e-12\\
\hline
$8\Delta t $&3.667e-05&9.2733e-10&4.6629e-12\\
\hline
$16\Delta t $&7.334e-05&1.8547e-09&9.3258e-12\\
\hline
\end{tabular}
\caption{\label{table_conv_7} Convergence tableau for the ABA-Splitting method.}

\end{table}
\begin{table}
\centering
\begin{tabular}{ |c|c|c|c| } 
 \hline
  & $\Delta x/4$ & $\Delta x/8$ & $\Delta x/16$  \\ 
\hline
$4\Delta t $&1.583e-05&4.2155e-10&2.2106e-12\\
\hline
$8\Delta t $&3.667e-05&9.2733e-10&4.6629e-12\\
\hline
$16\Delta t $&7.334e-05&1.8547e-09&9.3258e-12\\
\hline
\end{tabular}
\caption{\label{table_conv_8} Convergence tableau for the BAB-Splitting method. }
\end{table}

\begin{table}
\centering
\begin{tabular}{ |c|c|c|c| } 
 \hline
  & $\Delta x/4$ & $\Delta x/8$ & $\Delta x/16$  \\ 
\hline
$4\Delta t $&1.312e-05&3.3337e-10&1.6673e-12\\
\hline
$8\Delta t $&3.667e-05&9.2733e-10&4.6629e-12\\
\hline
$16\Delta t $&7.334e-05&1.8547e-09&9.3258e-12\\
\hline
\end{tabular}
\caption{\label{table_conv_9} Convergence tableau for the ABA(CN)Splitting method. }
\end{table}

\begin{table}
\centering
\begin{tabular}{ |c|c|c|c| } 
 \hline
 $$&$\Delta x/4$ & $\Delta x/8$ & $\Delta x/16$  \\ 
\hline
$4\Delta t $&1.312e-05&3.3337e-10&1.6673e-12\\
\hline
$8\Delta t $&3.667e-05&9.2733e-10&4.6629e-12\\
\hline
$16\Delta t $&7.334e-05&1.8547e-09&9.3258e-12\\
\hline
\end{tabular}
\caption{\label{table_conv_10} Convergence tableau for the ABA(semiCN) Splitting method. }
\end{table}

The computational times and the errors of the different methods
for the single soliton solutions are given in Table \ref{table_1} and \ref{table_2}.

\begin{table}
\centering
\begin{tabular}{ |c|c|c|c|c| } 
 \hline
 & T=2.5 & T=5 & T=7.5 & T=10 \\ 
\hline
 & T=2.5 & T=5 & T=7.5 & T=10 \\ 
 \hline
 Implicit Euler method & 0.8313& 1.6785 &2.1124 & 2.9281 \\ 
 \hline
 Crank-Nicolson scheme &2.0496 & 3.8930 & 5.7764 & 7.1148 \\ 
 \hline
 AB-splitting: A and B operators are spectral & 0.0271 & 0.0486 &0.0785 &0.1007 \\ 
 \hline
 AB-splitting: A Spectral , B FD & 1.8159 &3.2140 & 4.6932 & 5.7207  \\ 
 \hline
 AB-splitting: A FD , B Spectral & 0.0466 &  0.0551 & 0.0668& 0.0962 \\ 
 \hline
AB-splitting: A FD , B FD & 2.4798 & 3.8211 &5.7146 & 7.0136\\ 
 \hline
 ABA-Splitting & 0.0352&0.0632 & 0.0940 & 0.1264 \\ 
 \hline
 BAB-Splitting & 0.0343 &0.0624 &0.1003& 0.1281 \\ 
 \hline
 ABA(CN)-Splitting & 0.9762 &1.9774&2.6190&3.2304 \\ 
 \hline
 ABA(semiCN)-Splitting & 2.5906 &4.5933&6.5765&8.5612 \\ 
 \hline
\end{tabular}
\caption{\label{table_1} Computational times of one soliton with the different methods.}
\end{table}

\begin{table}
\centering
\begin{tabular}{ |c|c|c|c|c| } 
 \hline
 & T=2.5 & T=5 & T=7.5 & T=10 \\ 
  \hline
 Implicit Euler method & 0.8977& 1.9084&2.4616 & 2.6552 \\ 
 \hline
 Crank-Nicolson scheme&0.9165 & 2.0208& 2.7069 & 2.9975 \\ 
 \hline
 AB-splitting: A and B operators are spectral  &0.0330 & 0.0396& 0.0420 & 0.0488\\ 
 \hline
 AB-splitting: A Spectral , B FD & 0.9443 & 2.0468& 2.7054 &  2.9648 \\ 
 \hline
 AB-splitting: A FD , B Spectral & 0.1333 & 0.3893 & 0.7650 & 1.2144 \\ 
 \hline
AB-splitting: A FD , B FD & 0.9165 & 2.0208& 2.7069 & 2.9975\\ 
 \hline
 ABA-splitting &0.0057 &0.0080 & 0.0097 &0.0111 \\ 
 \hline
 BAB-splitting & 0.0057&0.0080 & 0.0097 & 0.0111 \\ 
 \hline
 ABA(CN)-Splitting &0.9178 &2.0201 &2.6952& 2.9630 \\ 
 \hline
 ABA(semiCN)-Splitting & 0.9174 &2.0208&2.7003&2.9740 \\ 
 \hline
\end{tabular}
\caption{\label{table_2} Numerical errors of one soliton with the different methods. }
\end{table}

\begin{figure}[ht]
\begin{center}  
\includegraphics[width=5.0cm,angle=-0]{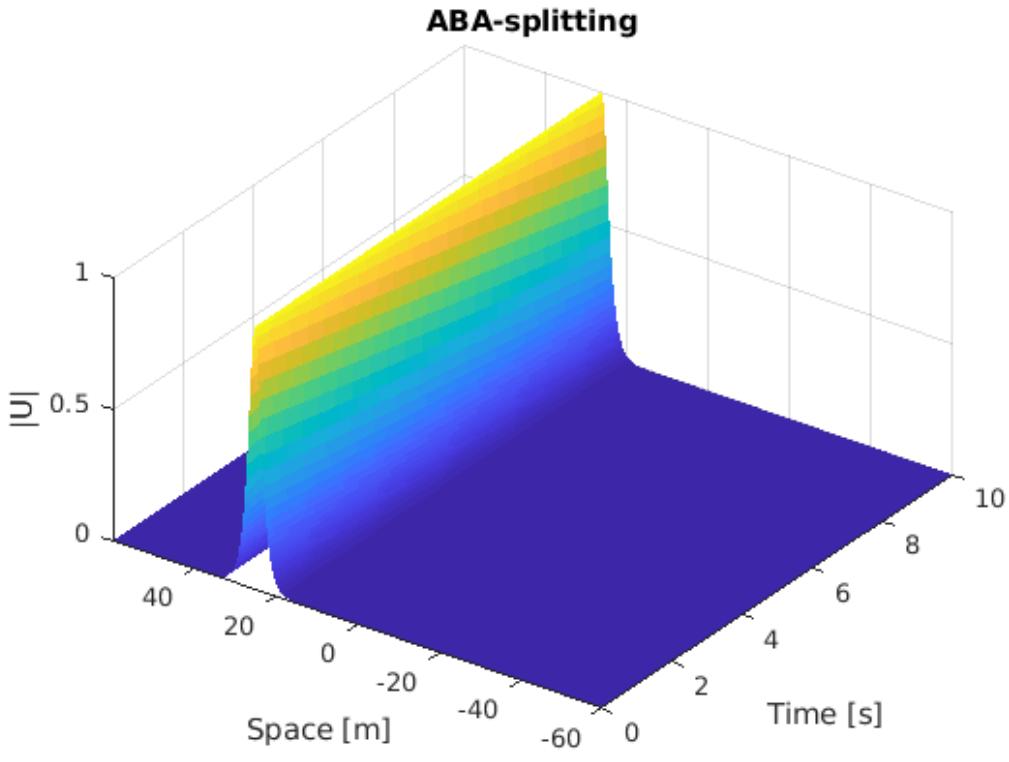}
\includegraphics[width=5.0cm,angle=-0]{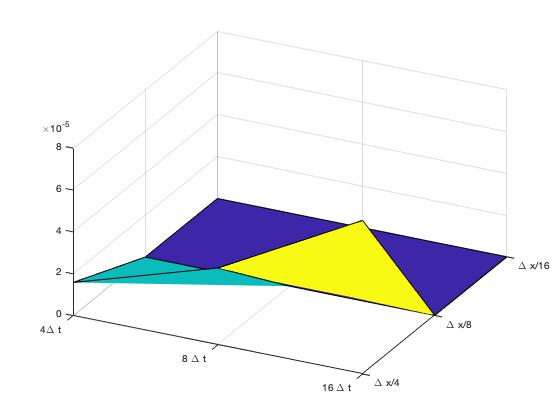}
\end{center}
\caption{\label{experiment_one}
Results of the GPE with one soliton equation, here we have appliedthe ABA-splitting approach (left figure: numerical results, right figure: convergence results).}
\end{figure}

The Figure \ref{experiment_one} present the solutions of the one soliton results and the convergence tableau.

\begin{figure}[ht]
\begin{center}  
\includegraphics[width=5.0cm,angle=-0]{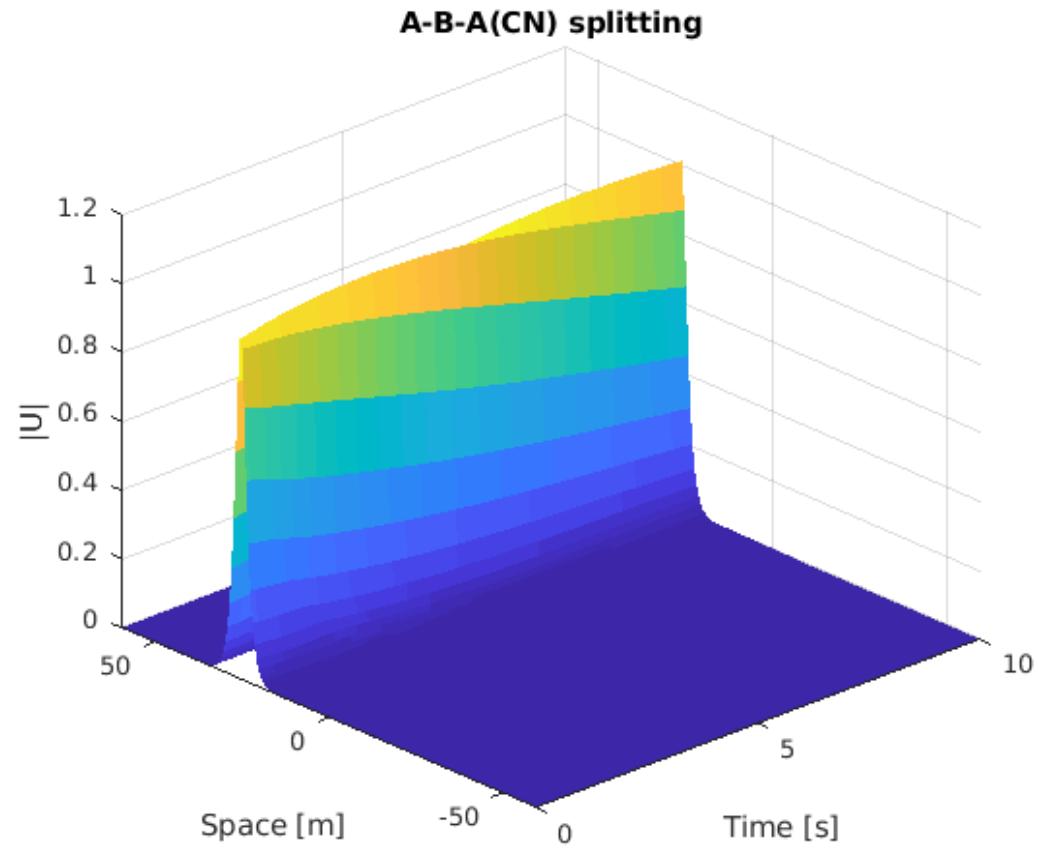} 
\end{center}
\caption{\label{experiment_one_2} Numerical solution with the ABA-CN method of the single solitons.}
\end{figure}

The Figure \ref{experiment_one_2} present the solutions with the approximated conservation finite difference scheme.

\begin{remark}
We see the benefits of the conservation schemes in the long time behaviour.
But the drawbacks are the time-consuming computations. The balance based on the
splitting approach including the conservative schemes are an alternative to
reduce the time-consuming approaches and allow to obtain asymptotic 
conservative results with sufficient enough iterative steps.
\end{remark}

\subsection{Second Example: Collision of two solitons}

We apply a collision of two solitons with the GPE.
The evolution equation is given as:
\begin{eqnarray}
&& \frac{\partial u}{\partial t} = - i H u , \; x \in \Omega, \; t \in [0, 10] , \\
&& u(x,0) = \sech(\frac{1}{\sqrt{2}} (x - 20)) \exp(-i \frac{x}{20}) + \\
&& +  \sech(x + 20) \exp(i \frac{x}{20}) , \; x \in \Omega ,  \\
&& u(x, t)= 0.0 , \; x \in \partial \Omega, \; t \in [0, T] ,
\end{eqnarray}
with $ H u = \left(-\frac{1}{2} \frac{\partial^2}{\partial x^2}  + g |u|^{2 \sigma} \right) u$, $\sigma = 1.0$.

We have two solitons starting in $x=-20$ and $x=20$ and they collide at $x = 0$ at the time-point $t=5.0$.

For the reference solution, we apply a fine spatial- and time-discretised solution with
an ABA method.

Further, we also decouple the full equation after the spatial discretisation into a linear and nonlinear operator part,
given as:
\begin{eqnarray}
&& H U^n = A(t,x, U^n) = A_1(t, x) + A_2(t,x, U^n) ,
\end{eqnarray}

In the Table \ref{table_3} and \ref{table_4}, we present the computational time and the
numerical errors of the different methods for the two-solitons modelling problem.
\begin{table}
\centering
\begin{tabular}{ |c|c|c|c|c| } 
 \hline
 & T=2.5 & T=5 & T=7.5 & T=10 \\ 
 \hline
 Implicit Euler method& 2.4928& 3.5601 &4.9031 & 6.2648 \\ 
 \hline
 Crank-Nicolson scheme &4.5648& 8.8923 & 13.8926 & 15.9353\\ 
 \hline
 AB-splitting: A and B operators are spectral  &0.0342 &0.0632 & 0.1004 &0.1429\\ 
 \hline
 AB-splitting: A Spectral , B FD & 3.5292 & 6.7374 & 9.7683&12.9375  \\ 
 \hline
 AB-splitting: A FD , B Spectral& 0.0349& 0.0678 & 0.0965& 0.1380\\ 
 \hline
AB-splitting: A FD , B FD &4.4182& 8.5995&12.3086& 16.4472\\ 
 \hline
 ABA-splitting & 0.0445&0.0858 & 0.1408 & 0.1989 \\ 
 \hline
 BAB-splitting & 0.0425 &0.0789 & 0.1524 & 0.1931 \\ 
  \hline
 ABA(CN)-Splitting & 2.1821 &4.4567 &6.3876& 7.7092 \\ 
 \hline
 ABA(semiCN)-Splitting & 6.1543&10.5217&16.1007&19.6879 \\ 
 \hline
\end{tabular}
\caption{\label{table_3} Computational times of two solitons with the different methods.}
\end{table}

\begin{table}
\centering
\begin{tabular}{ |c|c|c|c|c| } 
 \hline
 & T=2.5 & T=5 & T=7.5 & T=10 \\ 
 \hline
 Implicit Euler method& 1.0605& 4.6478 &5.0486 & 5.1546 \\ 
 \hline
 Crank-Nicolson scheme &1.0548& 4.3745& 9.1666 & 19.2207 \\ 
 \hline
 AB-splitting: A and B operators are spectral &0.0866 & 0.1059&0.1501 &0.1754 \\ 
 \hline
 AB-splitting: A Spectral , B FD &0.7579 &1.8421 & 2.4654 & 2.8003  \\ 
 \hline
 AB-splitting: A FD , B Spectral& 1.1001 & 6.0412 & 49.3173& 114.0526 \\ 
 \hline
AB-splitting: A FD , B FD &1.0548 & 4.3745& 9.1666& 19.2207\\ 
 \hline
 ABA-splitting & 0.0296 &0.0320& 0.0410& 0.0453 \\ 
 \hline
 BAB-splitting & 0.0295 &0.0314 &0.0405&0.0447 \\ 
  \hline
 ABA(CN)-Splitting &0.7024 &1.8489 &2.4949& 2.7952 \\ 
 \hline
 ABA(semiCN)-Splitting & 0.8599&2.6645&2.7771&2.7894\\ 
 \hline
\end{tabular}
\caption{\label{table_4} Numerical errors of two solitons with the different methods.}
\end{table}

\begin{figure}[ht]
\begin{center}  
\includegraphics[width=5.0cm,angle=-0]{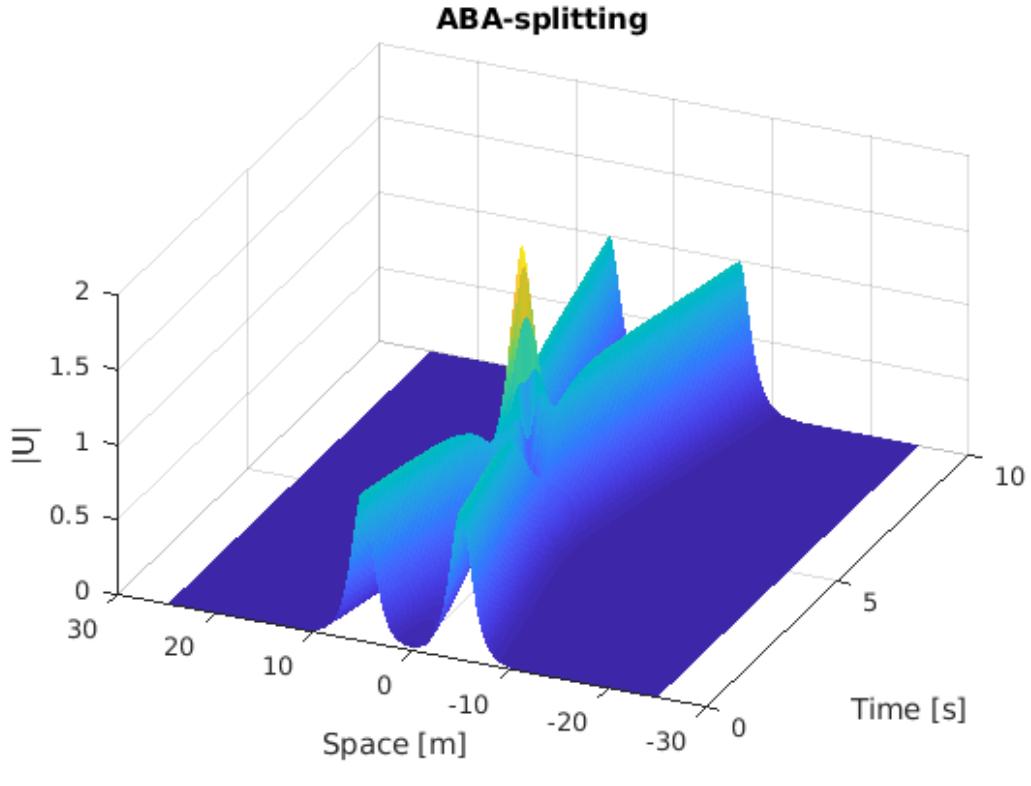} 
\includegraphics[width=5.0cm,angle=-0]{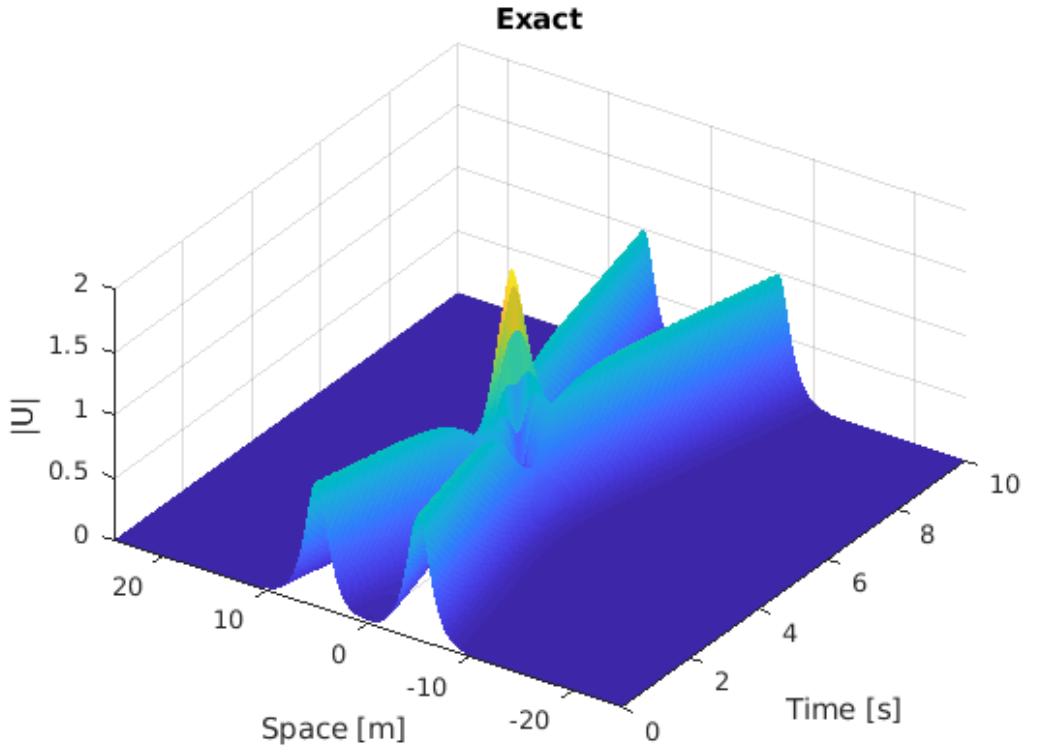}
\end{center}
\caption{\label{experiment_two}
 Results of the deterministic nonlinear Schr\"odinger equation with collisions of solitons (left figure: numerical results, right figure: exact results).}
\end{figure}

The Figure \ref{experiment_two} present the solutions and errors of the one soliton results.

\begin{figure}[ht]
\begin{center}  
\includegraphics[width=5.0cm,angle=-0]{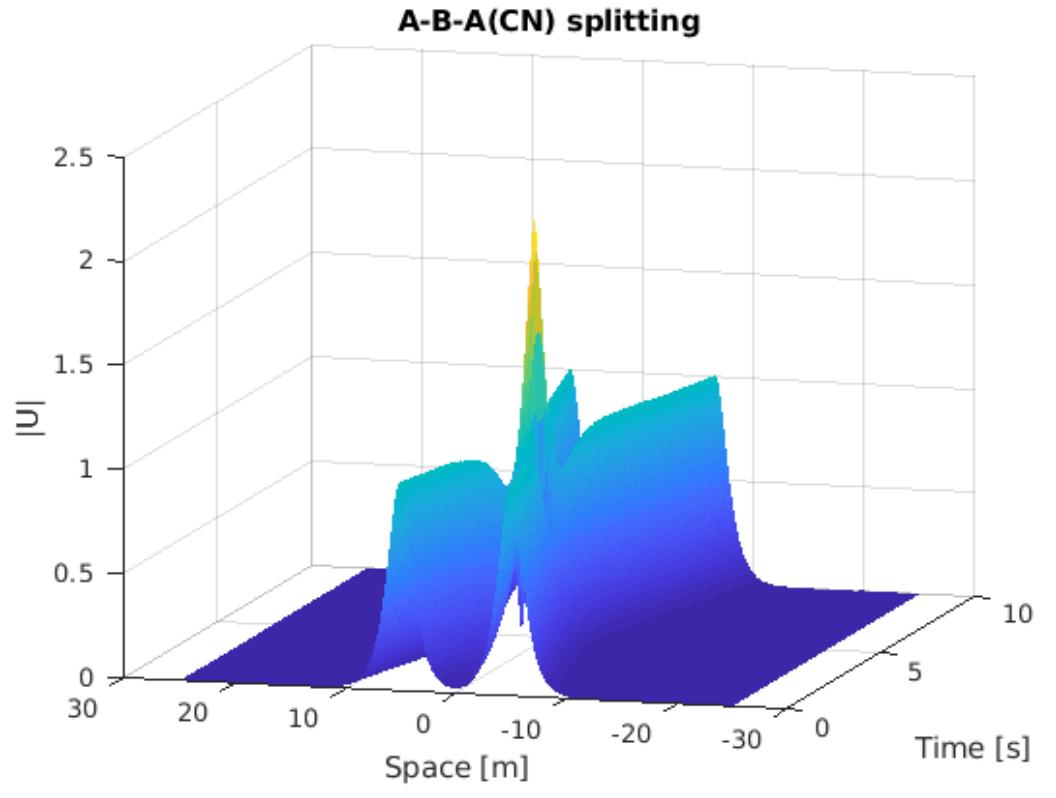} 
\end{center}
\caption{\label{experiment_4} Solution of the ABA-CN method for the two solitons.}
\end{figure}

The solution of the two-solitons with the ABA-CN method in Figure \ref{experiment_4}.

\begin{remark}
We also obtain the same results as for the single soliton solutions.
The alternative methods with the combination of the conservative schemes and the splitting approaches
have small numerical errors and optimal computational times in the area of the fast splitting methods.
With additional iterative steps, we could couple the ABA-iCN method more and achieve asymptotically
the conservation schemes.
\end{remark}

\section{Conclusion}
\label{concl}

We propose an alternative ABA-iCN method, which combines the conservative finite difference scheme with a fast ABA splitting approaches.
Such alternative methods allow to accelerate the solvers and stabilise the schemes to 
asymptotic conservative finite difference schemes.
We apply different numerical test examples and verify our assumptions. In future, we have to analyse carefully the
structure of the proposed methods with the underlying error analysis and present more
real-life applications in the field of soliton collisions.

\bibliographystyle{plain}

\end{document}